\documentclass[review]{elsarticle}

\usepackage{multirow}
\usepackage{amsmath,xfrac,cases,mathrsfs} 
\usepackage{amssymb}
\usepackage{enumerate}
\usepackage{lineno,hyperref}
\usepackage{braket}

\newtheorem{thm}{Theorem}
\newtheorem{lem}[thm]{Lemma}
\newdefinition{rmk}{Remark}
\newdefinition{prop}[thm]{Proposition}
\newproof{pf}{Proof}
\newproof{pot}{Proof of Theorem \ref{thm2}}

\journal{1}

\begin{document}
	
	\begin{frontmatter}
		
		\title{Drift parameter estimation for nonlinear reflected stochastic differential equations}

		\author[rvt]{Han Yuecai}
		\ead{hanyc@jlu.edu.cn}
		
		\author[rvt]{Zhang Dingwen\corref{cor1}}
		\ead{zhangdw20@mails.jlu.edu.cn}

		\cortext[cor1]{Corresponding author}

		\address[rvt]{School of Mathematics, Jilin University, Changchun, 130000, China}
		
		\begin{abstract}
		We study the maximum likehood estimator and least squares estimator for drift parameters of nonlinear reflected stochastic differential equations based on continuous observations. Under some regular conditions, we obtain the consistency and establish the asymptotic distributions of the two estimators. We briefly remark that our methods could be applied the the reflected stochastic processes with only one-sided reflecting barrier spontaneously.  Numerical studies show that the proposed estimators are adequate for practical use.
		\end{abstract}
		
		
	\end{frontmatter}
	
	\section{Introduction}\label{sec1}
	Let $(\Omega, \mathcal{F}, \mathbb{P}, \{\mathcal{F}_{t}\}_{t\geq0})$ be a filtered probability space and $W=\{W_{t}\}_{t\geq0}$ is a one-dimensional standard Brownian motion adapted to $\{\mathcal{F}_{t}\}_{t\geq0}$. The stochastic process $X=\{X_{t}\}_{t\geq0}$  is defined as the unique strong solution to the following nonlinear reflected stochastic differential equation (SDE) with two-sided barriers $l$ and $u$
	\begin{equation}\label{eqmodelT}
		\left\{
		\begin{aligned}
			&\mathrm{d}X_{t}=\theta f(X_{t})\mathrm{d}t+\sigma\mathrm{d}W_{t}+\mathrm{d}L_{t}-\mathrm{d}R_{t},\\
			&X_{0}=x\in[l,u],
		\end{aligned}
		\right.
	\end{equation}
	where $\sigma>0$, $0\leq l<u<\infty$ are given numbers, $f:[l,u]\rightarrow\mathbb{R}$ is a known measurable function and $\theta\in\mathbb{R}$ is a unknown parameter. The solution of a reflected SDE behaves like the solution of a standard SDE in the interior of its domain $(l,u)$. The processes $L=\{L_{t}\}_{t\geq0}$ and $R=\{R_{t}\}_{t\geq0}$ are the nimimal continuous increasing processes such that $X_{t}\in[l,u]$ for all $t\geq0$. Moreover, the two processes $L$ and $R$ subject to $L_{0}=R_{0}=0$ increase only when $X$ hits its boundary $l$ and $u$, and
	\begin{equation*}
		\int_{0}^{\infty}I_{\{X_{t}>l\}}\mathrm{d}L_{t}=0, \quad \int_{0}^{\infty}I_{\{X_{t}<u\}}\mathrm{d}R_{t}=0,
	\end{equation*} 
	where $I(\cdot)$ is the indicator function. For more about reflected stochastic processes, one can refer to \cite{Harrison2013}.
	
	Reflected SDEs have been widely used in many fields such as the queueing system \citep{Ward2003a,Ward2003b,Ward2005}, financial engineering \citep{Bo2010,Bo2011a,Han2016} and mathematical biology \citep{Ricciardi1987}. The reflected barrier is usually bigger than $0$ due to the physical restriction of the state processes which take non-negative values. One can refer to \cite{Harrison1985,Whitt2002} for more details on reflected SDEs and their broad applications. 
	
	The parameter estimation problem in reflected SDEs has gained much attention in recent years. A drift parameter estimator for a reflected fractional Brownian motion is proposed in \cite{Hu2013}. A maximum likelihood estimator (MLE) for the drift parameter based on the continuously observed reflected Ornstein–Uhlenbeck processes is proposed in \cite{Bo2011b}. A sequential MLE for the drift parameter of the reflected Ornstein–Uhlenbeck processes based on the observations throughout a random time interval is proposed in \cite{Lee2012}. The main tool of MLE is the Girsanov theorem for reflected stochastic processes. An ergodic type estimator for drift parameters based on discretely observed reflected O-U processes is proposed in \cite{Hu2015}. Subsequently, an ergodic type estimator for all parameters (drift and diffusion parameters) is proposed in \cite{Hu2021}. However, there is only limited literature on the drift parameter estimation of nonlinear reflected SDEs.

	The remainder of this paper is organized as follows. In section \ref{sec2}, we describe some notations and assumptions related to our context. In Section \ref{sec3}, we obtain the consistency and establish the asymptotic distributions of the two estimators. We also discuss the extension of main results to one-sided reflected processes. The proofs of main results are showed in Section \ref{sec4}. In section \ref{sec5}, we present some examples and give
	their numerical results. Section \ref{sec6} concludes with some discussion and remarks on the further work.

	\section{Preliminaries}\label{sec2}
	In this section, we present the necessary conditions and the ergodic properties for reflected stochastic processes. Throughout the paper, we shall use notation ``$\stackrel{P}{\longrightarrow}$" to denote ``convergence in probability" and notation ``$\sim$" to denote ``convergence in distribution”. 
	
	Now we introduce the following set of assumptions.
	\begin{enumerate}[i.]
		\item\label{asp1} The drift function $f(\cdot)$ satisfies a global Lipschitz condition, i.e., there exists a positive constant $M>0$ such that
		\begin{equation*}
			|f(x_{1})-f(x_{2})|\leq M|x_{1}-x_{2}|, \quad x_{1}, x_{2}\in [l,u].
		\end{equation*}
		\item\label{asp2}$f(\cdot)$ is not identically $0$ on $[l,u]$.
		\item\label{aspl} Let $\Theta=\{\theta:|\theta-\theta_{0}|\leq \rho\}$ be the parameter space for some constant $\rho$, where $\theta_{0}$ is the true parameter.
	\end{enumerate}
	
	Under condition \ref{asp1}, there exists a unique solution of Eq. (\ref{eqmodelT}) by an extension results of \cite{Lions1984}. 
	
	To this end, we give the unique invariant probability measure and the ergodic property for $X$.
	\begin{lem}\label{lemE}
		Under conditions \ref{asp1}-\ref{aspl}, for any function $b\in L_{1}([l,u])$, we have,
		\begin{enumerate}[a.]
			\item $X$ has a unique invariant probability measure,
			\begin{equation*}
				\pi_{l,u}(x)=\frac{e^{-\frac{2\theta}{\sigma^{2}} \int_{l}^{x}  f(y) \mathrm{d} y}}{\int_{l}^{u} e^{-\frac{2\theta}{\sigma^{2}} \int_{l}^{x}  f(y) \mathrm{d} y} \mathrm{d} x}.
			\end{equation*}
			\item The continuously observed process $\{X_{t}\}_{t\geq0}$ is ergodic,
			\begin{equation*}
				\lim_{T\rightarrow\infty}\frac{1}{T}\int_{0}^{T}b(X_{t})\mathrm{d}t=\mathbb{E}[b(X_{\infty})]=\int_{l}^{u}b(x)\pi_{l,u}(x)\mathrm{d}x.
			\end{equation*}
		\end{enumerate}
	\end{lem}
	With the unique invariant probability measure, we define 
	\begin{equation}\label{eqF}
		F=\int_{l}^{u}f^{2}(x)\pi_{l,u}(x)\mathrm{d}x.
	\end{equation}
	
	\section{Main Results}\label{sec3}
	\subsection{Maximum likehood estimator}
    In this subsection, we aim to construct the MLE for $\theta$ based on the continuous observations $\{X_{t}\}_{t\geq 0}$ and study its asymptotic behavior.
    
    The following lemma is the Girsanov’s formula for reflected Brownian motions.
    \begin{lem}\label{lemGir}
    Under the probability $\mathbb{P}^{0}$, $Y$ is a reflected Brownian motion satifies
    \begin{equation*}
    Y_{t}=\sigma\mathrm{d}W_{t}+\mathrm{d}L_{t}-\mathrm{d}R_{t},\quad Y_{0}=y \in[l,u].
    \end{equation*}
    Let 
    \begin{equation*}
    M(t,\gamma)=\exp\left(-\frac{\gamma}{\sigma}\int_{0}^{t}f(Y_{s})\mathrm{d}W_{s}-\frac{\gamma^{2}}{2\sigma^{2}}\int_{0}^{t}f^{2}(Y_{s})\mathrm{d}s\right).
    \end{equation*}
    Then, 
    \begin{enumerate}[a.]
    	\item For each $\gamma\in\mathbb{R}$, $\{M(t,\gamma)\}_{t\geq0}$ is a martingale adapted to the filtration generated by $W$.
    	\item Let $\mathbb{P}_{t}^{\gamma}=M(t,\theta)\mathrm{d}\mathbb{P}^{0}$, then $\mathbb{P}^{\theta}_{t}$ is a probability measure. Moreover, under $\mathbb{P}_{t}$, $Y$ is a reflected stochastic process with the drift term $\gamma f(y)$, diffusion term $\sigma$ and initial contion $y$. 
    \end{enumerate}
    \end{lem}
    
    In order to construct the MLE for $\theta_{0}$, we have the following null-hypotheses and alternative hypotheses
   \begin{equation*}
   \begin{aligned}
   &H_{0}: \text{The drift parameter is} \quad 0,\\
   &H_{1}: \text{The drift parameter is} \quad \theta\neq0.
   \end{aligned}
   \end{equation*} 
    Then the log-likehood function is given by 
    \begin{equation}\label{eqMLF}
    \begin{aligned}
    \ell_{T}(\theta)=&\frac{1}{T}\log\frac{\mathrm{d}\mathbb{P}^{\theta}}{\mathrm{d}\mathbb{P}^{0}}=\frac{1}{T}\left(\frac{\theta}{\sigma}\int_{0}^{T}f(X_{t})\mathrm{d}W_{t}+\frac{\theta^{2}}{2\sigma^{2}}\int_{0}^{T}f^{2}(X_{t})\mathrm{d}t\right)\\
    =&\frac{1}{T}\bigg(\frac{\theta}{\sigma^{2}}\int_{0}^{T}f(X_{t})\big(\mathrm{d}X_{t}-\mathrm{d}L_{t}+\mathrm{d}R_{t}\big)-\frac{\theta^{2}}{2\sigma^{2}}\int_{0}^{T}f^{2}(X_{t})\mathrm{d}t\bigg).
    \end{aligned}
    \end{equation}
    Then the MLE is defined as $\hat{\theta}_{T}=\arg\max_{\theta\in\Theta}\ell_{T}(\theta)$, which can be explicitly represent as 
    \begin{equation*}
    \hat{\theta}_{T}=\frac{\int_{0}^{T}f(X_{t})(\mathrm{d}X_{t}-\mathrm{d}L_{t}+\mathrm{d}R_{t})}{\int_{0}^{T}f^{2}(X_{t})\mathrm{d}t}.
    \end{equation*}

	The consistency of the MLE $\hat{\theta}_{n}$ is given as follows.
	\begin{thm}\label{thmConsMLE}
		Under conditions \ref{asp1}-\ref{aspl}, we have 
		\begin{equation*}
			\hat{\theta}_{T}(x){\longrightarrow} \theta_{0}
		\end{equation*}
		almost surely, as $T\rightarrow\infty$.
	\end{thm}

    The asymptotic distribution and rate of convergence of the MLE for the drift parameter is given as follows.
    \begin{thm}\label{thmAMLE}
    	Under conditions \ref{asp1}-\ref{aspl}, we have
    	\begin{equation*}
    		\sqrt{T}(\hat{\theta}_{T}-\theta_{0})\sim\mathcal{N}(0,\sigma^{2}F^{-1}),
    	\end{equation*}
    	as $T\rightarrow\infty$.
    \end{thm}

	\subsection{Least squares estimator}
	In this subsection, we aim to construct the LSE for $\theta$ based on the continuous obervations $\{X_{t}\}_{t\geq 0}$ and study its asymptotic behavior.
	
    This estimator is motivated by the following heuristic argument (\cite{Hu2010}). The LSE aims to minimize
    \begin{equation*}
    \int_{0}^{T}|\dot{X}_{t}-\theta f(X_{t})-\dot{L}_{t}+\dot{R}_{t}|^{2}\mathrm{d}t,
    \end{equation*}
    which is a quadratic function of $\theta$, though $\int_{0}^{T}(\dot{X}_{t}-\dot{L}_{t}+\dot{R}_{t})\mathrm{d}t$ does not exist. The minimum is achieved when 
    \begin{equation*}
    \tilde{\theta}_{T}=\frac{\int_{0}^{T}f(X_{t})(\mathrm{d}X_{t}-\mathrm{d}L_{t}+\mathrm{d}R_{t})}{\int_{0}^{T}f^{2}(X_{t})\mathrm{d}t}.
    \end{equation*}

	The consistency of the LSE $\tilde{\theta}_{T}$ is given as follows.
	\begin{thm}\label{thmConsLSE}
		Under conditions \ref{asp1}-\ref{aspl}, we have 
		\begin{equation*}
			\tilde{\theta}_{T}\stackrel{P}{\longrightarrow} \theta_{0},
		\end{equation*}
		as $n\rightarrow\infty$.
	\end{thm}
	
	The asymptotic distribution and rate of convergence of the LSE for the drift parameter is given as follows.
	\begin{thm}\label{thmALSE}
		Under conditions \ref{asp1}-\ref{aspl}, we have
		\begin{equation*}
		\sqrt{T}(\tilde{\theta}_{T}-\theta_{0})\sim\mathcal{N}(0,\sigma^{2}F^{-1}),
		\end{equation*}
		as $T\rightarrow\infty$.
	\end{thm}

	\subsection{Results of an expansion}
	In this subsection, we discuss the extension results of our main results to the reflected SDE with only one-sided barrier. Note that it is almost the same for only one reflecting barrier. We consider the following reflected SDE with a lower reflecting barrier $l$
	\begin{equation}\label{eqmodelo}
		\left\{
		\begin{aligned}
			&\mathrm{d}X_{t}=\theta f(X_{t})\mathrm{d}t+\sigma\mathrm{d}W_{t}+\mathrm{d}L_{t},\\
			&X_{0}=x\in[l,\infty).
		\end{aligned}
		\right.
	\end{equation}
	\begin{rmk}
		Our method can be applied to the reflected processes with only one-sided barrier Eq. (\ref{eqmodelo}). The unique invariant density of $X$ is given by
		\begin{equation*}
			\pi_{l}(x)=\frac{e^{-\frac{2}{\sigma^{2}} \int_{l}^{x} \theta f(y) \mathrm{d} y}}{\int_{l}^{\infty} e^{-\frac{2}{\sigma^{2}} \int_{l}^{x} \theta f(y) \mathrm{d} y} \mathrm{d} x}.
		\end{equation*}
		With the unique invariant density, we could do similar proofs as the proofs of Theorem \ref{thmConsMLE}, \ref{thmConsLSE}, \ref{thmAMLE} and \ref{thmALSE} to establish their consistency and asymptotic distribution. We omit the details here.
	\end{rmk}

	\section{Proofs of the Main Results}\label{sec4}
	In this section, we give proofs of the main results. We firstly give the proof of Lemma \ref{lemE} and \ref{lemGir}.
	
	\noindent\textbf{Proof of Lemma \ref{lemE}}. a. Note that $f(\cdot)$ is Lipschitz continuous and not identically $0$ on $[l,u]$. It suffices to verify that $X$ satisfies Conditions $5.2-5.5$ in \cite{Budhiraja2007}. Hence, $X$ has a unique invariant probability measure.  Define $m(x)$ and $s(x)$ the scale and speed densities
	\begin{equation*}
		s(x)=\exp \left(\frac{2\theta}{\sigma^{2}} \int_{l}^{x} f(r) \mathrm{d} r\right), \quad \text { and } \quad m(x)=\frac{2}{\sigma^{2} s(x)}.
	\end{equation*}
	If there exists a stationary density $\psi(x)$, it satisfies the Kolmogorov backward equation (\cite{Han2016,Karlin1981})
	\begin{equation*}
		\frac{\sigma^{2}}{2} \frac{\mathrm{d}^{2}}{\mathrm{d} x^{2}}(\psi(x))+\frac{\mathrm{d}}{\mathrm{d} x}(f(x) \psi(x))=0 .
	\end{equation*}
	Solving this equation yeilds 
	\begin{equation*}
		\psi(x) =m(x)\left(C_{1} \int_{l}^{x} s(y) \mathrm{d} y+C_{2}\right),
	\end{equation*}
	where we choose $C_{1}=0$ and $C_{2}=\big(\int_{l}^{u}m(x)\mathrm{d}x\big)^{-1}$. Then 
	\begin{equation*}
		\pi_{l,u}(x)=\frac{e^{-\frac{2\theta}{\sigma^{2}} \int_{l}^{x} f(y) \mathrm{d} y}}{\int_{l}^{u} e^{-\frac{2\theta}{\sigma^{2}} \int_{l}^{x} f(y) \mathrm{d} y} \mathrm{d} x}.
	\end{equation*}
	
	\noindent b. From that $X$ has a unique invariant probability measure, we can conclude that $\{X_{t},t\geq0\}$ is ergodic (\cite{Han2016,Budhiraja2007}).
	\hfill$\square$

	\noindent\textbf{Proof of Lemma \ref{lemGir}}. Let 
	\begin{equation*}
		G_{t}=-\frac{\theta}{\sigma}\int_{0}^{t}f(Y_{s})\mathrm{d}W_{s}-\frac{\theta^{2}}{2\sigma^{2}}\int_{0}^{t}f^{2}(Y_{s})\mathrm{d}s. 
	\end{equation*}
	It$\hat{o}$'s formula for reflected processes (see in \cite{Harrison1985}) yields
	\begin{equation*}
		e^{G_{t}}=1-\int_{0}^{t}\frac{\theta}{\sigma}f(Y_{s})\exp\bigg(-\frac{\theta}{\sigma}\int_{0}^{s}f(Y_{r})\mathrm{d}W_{r}-\frac{\theta^{2}}{2\sigma^{2}}\int_{0}^{r}f^{2}(Y_{r})\mathrm{d}r\bigg)\mathrm{d}W_{s}.
	\end{equation*}
	Let $\tau_{n}=\inf\{s\geq0:|\exp(G_{s})|\geq n\}$, it is obviously that $\tau_{n}$ is a stopping time and $\exp(G_{\tau_{n}\bigwedge s})\leq n$. Under conditions \ref{asp1}, we have that $f(X_{\tau_{n}})$ is bounded. Hence 
	$G_{\tau_{n}}$ is a martingale. Moreover, $\lim_{n\rightarrow\infty}\tau_{n}=\infty$ means that $G_{t}$ is a local martingale. Moreover, $\mathbb{E}(G_{t})=1$, which completes the proof.
	\hfill$\square$
	
	Before we give the proofs of the main results, we prepare some preliminary lemmas. Recall the $F$ defined in Eq. (\ref{eqF}).
	
	\begin{lem}\label{lem1}
	Under conditions \ref{asp1}-\ref{aspl}, we have
	\begin{equation*}
	\frac{1}{T}\int_{0}^{T}f^{2}(X_{t})\mathrm{d}t\longrightarrow F
	\end{equation*}
    almost surely, as $T\rightarrow\infty$.
	\end{lem}
	\noindent\textbf{Proof of Lemma \ref{lem1}}. The lemma is an immediate consequence of Lemma \ref{lemE}.
	
	\begin{lem}\label{lem2}
	Under conditions \ref{asp1}-\ref{aspl}, we have
	\begin{equation*}
	\frac{1}{T}\int_{0}^{T}f(X_{t})\mathrm{d}W_{t}\stackrel{P}{\longrightarrow}0,
	\end{equation*}
    as $T\rightarrow\infty$.
	\end{lem}
    \noindent\textbf{Proof of Lemma \ref{lem2}}. Note that 
    \begin{equation*}
    \mathbb{E}\big[\frac{1}{T}\int_{0}^{T}f(X_{t})\mathrm{d}W_{t}\big]=0.
    \end{equation*}
	By It$\hat{o}$ isometry and the the dominated convergence theorem, we have
	\begin{equation*}
	\begin{aligned}
	\lim_{T\rightarrow\infty}\mathbb{E}\big[\frac{1}{T}\int_{0}^{T}f(X_{t})\mathrm{d}W_{t}\big]^{2}=&\lim_{T\rightarrow\infty}\frac{1}{T^{2}}\int_{0}^{T}\mathbb{E}[f^{2}(X_{t})]\mathrm{d}t\\
	=&\lim_{T\rightarrow\infty}\frac{1}{T^{2}}\int_{0}^{T}f^{2}(X_{t})\mathrm{d}t.
	\end{aligned}
	\end{equation*}
	By Lemma \ref{lem1}, we have 
	\begin{equation*}
	\begin{aligned}
	\lim_{T\rightarrow\infty}\frac{1}{T^{2}}\int_{0}^{T}f^{2}(X_{t})\mathrm{d}t=\lim_{T\rightarrow\infty}\frac{F}{T}=O(T^{-1}).
	\end{aligned}
	\end{equation*}
	By Chebyshev's inequality, we obtain the desired results.
	\hfill$\square$

	\begin{lem}\label{lem3}
	 Let
	\begin{equation*}
	H_{T}=\frac{1}{\sqrt{T}}\int_{0}^{T}f(X_{t})\mathrm{d}W_{t}.
	\end{equation*}
    Under conditions \ref{asp1}-\ref{aspl}, we have
	\begin{equation*}
	H_{T}\sim \mathcal{N}(0, F),
	\end{equation*}
    as $T\rightarrow\infty$.
	\end{lem}
	\noindent\textbf{Proof of Lemma \ref{lem3}}. By Lemma \ref{lem1}, we have 
	\begin{equation*}
	\frac{1}{T}\int_{0}^{T}f^{2}(X_{t})\mathrm{d}t\rightarrow F.
	\end{equation*}
    It follows immediately from $X_{t}$ is adapted with respect to $\mathcal{F}_{t}$ that $H_{T}$ converges in law as $T\rightarrow\infty$ to a centered normal distribution with variance $F$.
    \hfill$\square$
	
	Before we give the proofs of Theorem \ref{thmConsMLE} and \ref{thmAMLE}, we do some strightforward caculation. By the log-likehood function Eq. (\ref{eqMLF}), we have
	\begin{equation*}
	\begin{aligned}
	\ell^{\prime}_{T}(\theta)=&\frac{1}{T}\bigg(\frac{1}{\sigma^{2}}\int_{0}^{T}f(X_{t})\big(\mathrm{d}X_{t}-\mathrm{d}L_{t}+\mathrm{d}R_{t}\big)-\frac{\theta}{\sigma^{2}}\int_{0}^{T}f^{2}(X_{t})\mathrm{d}t\bigg),\\
	\end{aligned}
	\end{equation*}
	and
	\begin{equation*}
	\ell_{T}^{\prime\prime}(\theta)=-\frac{1}{T\sigma^{2}}\int_{0}^{T}f^{2}(X_{t})\mathrm{d}t.
	\end{equation*}
	\noindent\textbf{Proof of Theorem \ref{thmConsMLE}}. Note that
	\begin{equation*}
	\begin{aligned}
	\ell_{T}(\theta)=&\frac{1}{T}\bigg(\frac{\theta}{\sigma^{2}}\int_{0}^{T}f(X_{t})\big(\theta_{0}f(X_{t})\mathrm{d}t+\sigma\mathrm{d}W_{t}\big)-\frac{\theta^{2}}{2\sigma^{2}}\int_{0}^{T}f^{2}(X_{t})\mathrm{d}t\bigg)\\
	=&\frac{1}{T}\bigg(\frac{2\theta\theta_{0}-\theta^{2}}{2\sigma^{2}}\int_{0}^{T}f^{2}(X_{t})\mathrm{d}t+\frac{\theta}{\sigma}\int_{0}^{T}f(X_{t})\mathrm{d}W_{t}\bigg).
	\end{aligned}
	\end{equation*}
    Let
    \begin{equation*}
    D(\theta)=\frac{2\theta\theta_{0}-\theta^{2}}{2\sigma^{2}}F.
    \end{equation*}
	By Lemma \ref{lem1} and \ref{lem2}, we have
	\begin{equation*}
	\ell_{T}(\theta)\stackrel{P}{\longrightarrow}D(\theta),
	\end{equation*}
	as $T\rightarrow\infty$. Moreover, $D(\theta)$ is uniquely maximized when $\theta=\theta_{0}$. 
	
	Similarly to the proof of Lemma $3$ in \cite{Amemiya1973}, we conclude that $\hat{\theta}_{T}\rightarrow\theta_{0}$ almost surely, as $T\rightarrow\infty$.
	\hfill$\square$

	\noindent\textbf{Proof of Theorem \ref{thmAMLE}}. Note that
	\begin{equation*}
	\ell_{T}^{\prime}(\hat{\theta}_{T})=\ell_{T}^{\prime}(\theta_{0})+\ell_{T}^{\prime\prime}(\bar{\theta})(\hat{\theta}_{T}-\theta_{0})+o\big((\hat{\theta}-\theta_{0})^{2}\big),
	\end{equation*}
	where $\bar{\theta}$ is a mean value, located between $\hat{\theta}_{T}$ and $\theta_{0}$. By the consistency of $\hat{\theta}_{T}$, we have $\bar{\theta}=\theta_{0}$. Moreover, $\ell_{T}^{\prime\prime}(\bar{\theta})=\ell_{T}^{\prime\prime}(\theta_{0})$. Then,
	\begin{equation*}
	\begin{aligned}
	\sqrt{T}\big(\hat{\theta}_{T}-\theta_{0}\big)=-\frac{\ell_{T}^{\prime}(\theta_{0})}{\ell_{T}^{\prime\prime}(\theta_{0})}=\frac{\frac{\sigma}{\sqrt{T}}\int_{0}^{T}f(X_{t})\mathrm{d}W_{t}}{\frac{1}{T}\int_{0}^{T}f^{2}(X_{t})\mathrm{d}t}.
	\end{aligned}
	\end{equation*}
	By Lemma \ref{lem3}, we have
	\begin{equation*}
	\frac{\sigma}{\sqrt{T}}\int_{0}^{T}f(X_{t})\mathrm{d}W_{t}\sim\mathcal{N}(0,\sigma^{2}F).
	\end{equation*}
    By Lemma \ref{lem1} and Slutsky's theorem, we obtain the desired results.
	\hfill$\square$

	Before we give the proofs of Theorem \ref{thmConsLSE} and \ref{thmALSE}, we propose a useful alternative expression for the LSE
	\begin{equation}\label{eqalter1}
		\tilde{\theta}_{T}=\frac{\int_{0}^{T}f(X_{t})(\mathrm{d}X_{t}-\mathrm{d}L_{t}+\mathrm{d}R_{t})}{\int_{0}^{T}f^{2}(X_{t})\mathrm{d}t}=\theta_{0}+\frac{\sigma\int_{0}^{T}f(X_{t})\mathrm{d}W_{t}}{\int_{0}^{T}f^{2}(X_{t})\mathrm{d}t}.
	\end{equation}
	\noindent\textbf{Proof of Theorem \ref{thmConsLSE}}.	
	By Lemma \ref{lem1} and \ref{lem2}, we have
	\begin{equation}\label{eqalter2}
	\frac{\frac{\sigma}{T}\int_{0}^{T}f(X_{t})\mathrm{d}W_{t}}{\frac{1}{T}\int_{0}^{T}f^{2}(X_{t})\mathrm{d}t}\stackrel{P}{\longrightarrow}0,
	\end{equation}
	as $T\rightarrow\infty$. Then we conclude the results combining Eq. (\ref{eqalter1}) and (\ref{eqalter2}).
	\hfill$\square$
	
	\noindent\textbf{Proof of Theorem \ref{thmALSE}}. By Eq. (\ref{eqalter1}), we have
	\begin{equation*}
	\sqrt{T}(\tilde{\theta}_{T}-\theta_{0})=\frac{\frac{\sigma}{\sqrt{T}}\int_{0}^{T}f(X_{t})\mathrm{d}W_{t}}{\frac{1}{T}\int_{0}^{T}f^{2}(X_{t})\mathrm{d}t}.
	\end{equation*}
	By Lemma \ref{lem3}, we have
	\begin{equation*}
	\frac{\sigma}{\sqrt{T}}\int_{0}^{T}f(X_{t})\mathrm{d}W_{t}\sim\mathcal{N}(0,\sigma^{2}F)
	\end{equation*}
    as $T\rightarrow\infty$. 
    
    By Lemma \ref{lem1} and Slutsky's theorem, we obtain the desired results.
    \hfill$\square$
	\begin{rmk}\label{rmk2}
	Note that the MLE and LSE have the same form through some transformations. The proofs of consistency and saymptotic distribution for LSE could be applied to the proofs of MLE. However, we have used traditional statistical approach here. 
	\end{rmk}
	\section{Numerical Results}\label{sec5}
	In this section, we present several examples on the cases of two-sided barriers and one-sided barrier for illustration. We take the lower reflecting barrier $l=0$ and the upper reflecting barrier $u=3$. 
	The diffusion parameter $\sigma=2$ and the drift function we choose shall satisfy the condition \ref{asp1} and \ref{asp2}.  Hence, we considered the following cases 
	\begin{enumerate}[a.]
		\item $f_{1}(x)=\sin(2\pi x)$
		\item $f_{2}(x)=\sqrt{x}$.
	\end{enumerate}
	
	We examine the numerical behavior of the proposed estimator for reflected SDEs (\ref{eqmodelT}) and (\ref{eqmodelo}). For a simulation of reflected stochastic processes, we make use of the numerical method presented in \cite{Lepingle1995}, which is known to yield the same rate of convergence as the usual Euler-Maruyama method.
	
	For each setting, we generate $N=1000$ Monte Carlo simulations of the sample paths, each consisting of $n=100, 200$ and $500$ observations. Further simulations evidence that the use of other drift function has little impact on the estimator's empirical performance, we ommit the results here. 
	The proposed estimator is evaluated by the Bias, the sample standard deviation (Std.dev) and the mean squares error (MSE).

	Table \ref{table1} and \ref{table2} summarize the main findings over 1000 simulations. We observe that as the sample size increases, the Bias decreases and is small, that the empirical and model-based standard errors agree reasonably. The performance improves with larger sample sizes.

	\begin{table}
		\caption{Simulation results with different sample size $n$.}
		\begin{tabular}{cccccccc}
			\hline 
			\multicolumn{8}{c}{Case 1: $f_{1}(x)=\sin(2\pi x)$.} \tabularnewline
			\hline 
			\hline 
			&&\multicolumn{3}{c}{Two-sided}
			&\multicolumn{3}{c}{One-sided} 
			\tabularnewline
			\cline{3-5}  \cline{6-8}
			\makebox[0.06\textwidth]{$\theta_{0}$}&
			\makebox[0.06\textwidth]{n} &
			\makebox[0.1\textwidth]{Bias}&  
			\makebox[0.1\textwidth]{Std.dev}&
			\makebox[0.1\textwidth]{MSE}&
			\makebox[0.1\textwidth]{Bias}&  
			\makebox[0.1\textwidth]{Std.dev}&
			\makebox[0.1\textwidth]{MSE}
			\tabularnewline
		     1 & 100 & 0.076 & 0.460 & 0.217
			   & 0.063  & 0.477 & 0.231\tabularnewline
			   & 200 & 0.049 & 0.400 & 0.163
			   & 0.057  & 0.315 & 0.103\tabularnewline
		       & 500 & 0.035 & 0.305 & 0.095
			   & 0.036  & 0.227 & 0.053\tabularnewline
		    2  & 100 & 0.107 & 0.476 & 0.238
			   & 0.102  & 0.491 & 0.251\tabularnewline
			   & 200 & 0.079 & 0.405 & 0.170
			   & 0.087  & 0.413 & 0.178\tabularnewline
			   & 500 & 0.036 & 0.293 & 0.087
			   & 0.039  & 0.286 & 0.084\tabularnewline
			\hline 
		\end{tabular}
		\label{table1}
	\end{table}
	
	\begin{table}
		\caption{Simulation results with different sample size $n$.}
		\begin{tabular}{cccccccc}
			\hline 
			\multicolumn{8}{c}{Case 2: $f_{2}(x)=\sqrt{x}$.} \tabularnewline
			\hline 
			\hline 
			&&\multicolumn{3}{c}{Two-sided}
			&\multicolumn{3}{c}{One-sided} 
			\tabularnewline
			\cline{3-5}  \cline{6-8}
			\makebox[0.06\textwidth]{$\theta_{0}$}&
			\makebox[0.06\textwidth]{n} &
			\makebox[0.1\textwidth]{Bias}&  
			\makebox[0.1\textwidth]{Std.dev}&
			\makebox[0.1\textwidth]{MSE}&
			\makebox[0.1\textwidth]{Bias}&  
			\makebox[0.1\textwidth]{Std.dev}&
			\makebox[0.1\textwidth]{MSE}
			\tabularnewline
			1 & 100  & -0.255 & 0.716 & 0.577
			         & -0.235  & 0.683 & 0.521\tabularnewline
			  & 200  & -0.066 & 0.283 & 0.085
			         & -0.077  & 0.271 & 0.079\tabularnewline
			  & 500  & -0.040 & 0.073 & 0.007
			         & -0.010  & 0.063 & 0.004\tabularnewline
			2 & 100  & -0.076 & 0.392 & 0.160
			         & -0.088  & 0.382 & 0.154\tabularnewline
			  & 200  & -0.041 & 0.124 & 0.017
			         & -0.013  & 0.125 & 0.016\tabularnewline
			  & 500  & -0.024 & 0.061 & 0.004
			         & -0.001  & 0.030 & 0.001\tabularnewline
			\hline 
		\end{tabular}
		\label{table2}
	\end{table}
	\section{Conclusion}\label{sec6}
	We provided the explicit MLE and LSE for the drift parameter of two-sided reflceted stochastic processes based on continuously observed processes. We obtain the consistency and establish the asymptotic distributions of the two estimators. The proposed estimators could be applied to reflected stochastic processes with only one-sided reflecting barrier. A simulation study is presented to show that our estimators are adequate for practical use.
	
	Some further researches may include providing the other estimators for the reflected stochastic processes and investigating the statistical inference for the other reflected diffusions.


\begin{thebibliography}{99}
		\footnotesize
		\bibitem{Amemiya1973}{Amemiya, T., 1973. Regression Analysis When the Dependent Variable is Truncated Normal. Econometrica 41, 997-1016.}
		
		\bibitem{Budhiraja2007}{Budhiraja, A., Lee, C., 2007. Long time asymptotics for constrained diffusions in polyhedral domains. Stoch. Process. Their Appl. 117(8), 1014–1036.}
		
		\bibitem[Bo et al. (2010)]{Bo2010}{Bo, L., Wang, Y., Yang, X., 2010. Some integral functionals of reflected SDEs and their applications in finance. Quant. Financ. 11, 343–348.}
		
		\bibitem[Bo et al. (2011a)]{Bo2011a}{Bo, L., Tang, D., Wang, Y., Yang, X., 2011a. On the conditional default probability in a regulated market: a structural approach. Quant. Financ. 11, 1695–1702.}
		
		\bibitem[Bo et al. (2011b)]{Bo2011b}{Bo, L., Wang, Y., Yang, X., Zhang, G., 2011b. Maximum likelihood estimation for reflected Ornstein–Uhlenbeck processes. J. Stat. Plan. Infer. 141, 588–596.}
		
		\bibitem{Frydman1980}{Frydman, R., 1980. A proof of the consistency of maximum likelihood estimators of non-linear regression models with autocorrelated errors. Econometrica 48, 853-860.}
		
		\bibitem[Harrison (1985)]{Harrison1985}{Harrison, J.M., 1985. Brownian motion and stochastic flow systems. Wiley, New York.}
		
		\bibitem{Hu2010}{Hu, Y., Nualart, D., 2010. Parameter estimation for fractional Ornstein–Uhlenbeck processes. Stat. Probab. Lett. 80, 1030-1038.}
		
		\bibitem{Harrison2013}{Harrison, J.M., 2013. Brownian Models of Performance and Control. Cambridge University Press, New York.}
		
		\bibitem{Hu2013}{Hu, Y., Lee, C., 2013. Drift parameter estimation for a reflected fractional Brownian motion based on its local time. J. Appl. Probab. 50, 592-597.}
		
		\bibitem[Hu et al. (2015)]{Hu2015}{Hu, Y., Lee, C., Lee, M. H., Song, J., 2015. Parameter estimation for reflected Ornstein-Uhlenbeck processes with discrete observations. Stat. Infer. Stoch. Proc. 18, 279–291.}
		
		\bibitem[Han et al. (2016)]{Han2016}{Han, Z, Hu, Y., Lee, C., 2016. Optimal pricing barriers in a regulated market using reflected diffusion processes. Quant. Financ. 16, 639-647.}
		
		\bibitem[Hu and Xi (2021)]{Hu2021}{Hu, Y., Xi, Y., 2021. Estimation of all parameters in the reflected Ornstein–Uhlenbeck process from discrete observations. Stat. Probab. Lett. 174.}
		
		\bibitem{Karlin1981}{Karlin, S., Taylor, H.M., 1981. A Second Course in Stochastic Processes, first ed. Academic Press, New York.}
		
		\bibitem{Kasonga1988}{Kasonga, R.A., 1988. The consistency of a non-linear least squares estimator from diffusion processes. Stoch. Process. Their Appl. 30, 263-275.}
		
		\bibitem{Lions1984}{Lions, P.L., Sznitman, A.S., 1984. Stochastic differential equations with reflecting boundary conditions. Commun. Pure Appl. Math. 37, 511-537.}
		
		\bibitem{Longstaff1989}{Longstaff, F.A., 1989. A nonlinear  general equilibrium model of the term structure of interest   rates. J. Financ. Econ. 23, 195-224.}
		
		\bibitem{Lepingle1995}{L\'epingle, D., 1995. Euler scheme for reflected stochastic differential equations. Math. Comput. Simul. 38, 119–126.}
		
		\bibitem[Lee et al. (2012)]{Lee2012}{Lee, C., Bishwal, J. P. N., Lee, M. H., 2012. Sequential maximum likelihood estimation for reflected Ornstein–Uhlenbeck processes. J. Stat. Plan. Infer. 142, 1234–1242.}
		
		\bibitem[Ricciardi and Sacerdote (1987)]{Ricciardi1987}{Ricciardi, L.M., Sacerdote, L., 1987. On the probability densities of an Ornstein–Uhlenbeck process with a reflecting boundary. J. Appl. Probab. 24, 355–369.}
		
		\bibitem[Whitt (2002)]{Whitt2002}{Whitt, W., (2002). Stochastic-process limits. Springer, New York.}
		
		\bibitem[Ward and Glynn (2003a)]{Ward2003a}{Ward, A. R., Glynn, P.W., 2003a. A diffusion approximation for a Markovian queue with reneging. Queueing Syst. 43, 103–128.}
		
		\bibitem[Ward and Glynn (2003b)]{Ward2003b}{Ward, A.R., Glynn, P.W., 2003b. Properties of the reflected Ornstein–Uhlenbeck process. Queueing Syst. 44, 109–123.}
		
		\bibitem[Ward and Glynn (2005)]{Ward2005}{Ward, A. R., Glynn, P. W., 2005. A diffusion approximation for a GI=GI=1 queue with balking or reneging. Queueing Syst. 50, 371–400.}
		
	\end{thebibliography}
\end{document}